\documentclass{gtart_h}

\def\ifplaintex{\expandafter\ifx\csname documentclass\endcsname\relax}

\def\gtp{{\mathsurround=0pt\it $\cal G\mskip-2mu$eometry \&\ 
$\cal T\!\!$opology $\cal P\!$ublications}}  

\def\recd{{\small Received:\qua\receiveddate\ifx\reviseddate\relax
\else\qquad Revised:\qua\reviseddate\fi\par}} 

\def\lognumber#1{\def\thelognumber{#1}}
\def\volumenumber#1{\def\thevolumenumber{#1}}
\def\volumeyear#1{\def\thevolumeyear{#1}}
\def\papernumber#1{\def\thepapernumber{#1}}
\def\pagenumbers#1#2{\def\startpage{#1}\def\finishpage{#2}}
\def\published#1{\def\publishdate{#1}}

\def\received#1{\def\receiveddate{#1}}

\def\accepted#1{\def\accepteddate{#1}}
\def\asciititle#1{\def\theasciititle{#1}}

\long\def\asciiabstract#1{\long\def\theasciiabstract{#1}}
\def\asciikeywords#1{\def\theasciikeywords{#1}}

\let\\\par\let\thelognumber\relax\let\thevolumenumber\relax
\let\thepapernumber\relax\let\thevolumeyear\relax\let\startpage\relax
\let\finishpage\relax\let\publishdate\relax\let\receiveddate\relax
\let\reviseddate\relax\let\accepteddate\relax\let\theasciititle\relax
\let\theasciiauthors\relax
\let\theasciiabstract\relax\let\theasciikeywords\relax

\let\theasciiemail\relax

\ifplaintex
\font\logobig=cmssbx10 scaled 3836
\font\logomed=cmssbx10 scaled 2557
\else
\font\logobig=cmssbx10 scaled 4200
\font\logomed=cmssbx10 scaled 2800
\fi

\long\def\makeagttitle{   
\count0=\startpage
\agt\hfill      
\hbox to 45truept{\vbox to 0pt{\vglue -13truept{\logomed A\kern -.37em{\logobig 
T}\kern -.38em G}\vss}\hss}
\break
{\small Volume \thevolumenumber\ (\thevolumeyear)
\startpage--\finishpage\nl
Published: \publishdate}

\vglue .25truein

{\parskip=0pt\leftskip 0pt plus
1fil\def\\{\par\smallskip}{\Large\bf\thetitle}\par\medskip} \vglue
0.05truein

{\parskip=0pt\leftskip 0pt plus 1fil\def\\{\par}{\sc\theauthors}
\par\medskip}%
 
\vglue 0.03truein 

{\small\leftskip 25truept\rightskip 25truept{\bf Abstract}\stdspace\theabstract

{\bf AMS Classification}\stdspace\theprimaryclass
\ifx\thesecondaryclass\relax\else; \thesecondaryclass\fi\par
{\bf Keywords}\stdspace \thekeywords\par}\vglue 7truept

}   

\ifplaintex
\font\phead=cmsl9 scaled 950
\font\pnum=cmbx10 scaled 913
\font\pfoot=cmsl9 scaled 950
\headline{\vbox to 0pt{\vskip -4.5mm\line{\small\phead\ifnum
\count0=\startpage ISSN 1472-2739 (on-line) 1472-2747 (printed)
\hfill {\pnum\folio}\else\ifodd\count0\def\\{ }%
\ifx\theshorttitle\relax\thetitle\else\theshorttitle\fi\hfill{\pnum\folio}
\else\def\\{ and }{\pnum\folio}\hfill\ifx\theshortauthors\relax\theauthors
\else\theshortauthors\fi\fi\fi}\vss}}
\footline{\vbox to 0pt{\vglue 0mm\line{\small\pfoot\ifnum\count0=\startpage
\copyright\ \gtp\hfill\else
\agt, Volume \thevolumenumber\ (\thevolumeyear)\hfill\fi}\vss}}
\else
\font=cmsl9 scaled 1050
\font\lnum=cmbx10 
\font=cmsl9 scaled 1050
\makeatletter
\def\@oddhead{{\small\lhead\ifnum\count0=\startpage ISSN 1472-2739 
(on-line) 1472-2747 (printed)\hfill {\lnum\number\count0}\else\ifodd\count0
\def\\{ }\ifx\theshorttitle\relax \thetitle \else\theshorttitle\fi\hfill
{\lnum\number\count0}\else\def\\{ and }{\lnum\number\count0}
\hfill\ifx\theshortauthors\relax 
\theauthors\else\theshortauthors\fi\fi\fi}}\def\@evenhead{\@oddhead}
\def\@oddfoot{\small\lfoot\ifnum\count0=\startpage\copyright\ \gtp\hfill\else
\agt, Volume \thevolumenumber\ (\thevolumeyear)\hfill\fi}
\def\@evenfoot{\@oddfoot}
\makeatother
\fi
\let\maketitlepage\makeagttitle
\let\makeshorttitle\maketitlepage
\let\maketitle\maketitlepage

\newwrite\gtoutfile
\long\gdef\makeheadfile{  
{\def\\{, }\def\s{ }
\immediate\openout\gtoutfile head.xxx
\immediate\write\gtoutfile{Proxy-for: \ifx\theasciiauthors\relax
\theauthors\else\theasciiauthors\fi\s<\ifx\theasciiemail\relax\theemail\else\theasciiemail\fi>}
\immediate\write\gtoutfile{\noexpand\\}
\immediate\write\gtoutfile{Authors: \ifx\theasciiauthors\relax
\theauthors\else\theasciiauthors\fi}
{\def\\{ }\immediate\write\gtoutfile{Title: \ifx\theasciititle\relax
\thetitle\else\theasciititle\fi}}
\immediate\write\gtoutfile{Subj-class: GT or SG, GR etc}
\immediate\write\gtoutfile{MSC-class: \theprimaryclass\ifx\thesecondaryclass\relax\else, \thesecondaryclass\fi}
\immediate\write\gtoutfile{Journal-ref: Algebr. Geom. Topol. \thevolumenumber\s
(\thevolumeyear) \startpage-\finishpage}
\immediate\write\gtoutfile{Comments: Published by Algebraic and
Geometric Topology at}
\immediate\write\gtoutfile{\s\s\s  http://www.maths.warwick.ac.uk/agt/AGTVol\thevolumenumber/agt-\thevolumenumber-\thepapernumber.abs.html}
\immediate\write\gtoutfile{\noexpand\\}
\immediate\write\gtoutfile{}
\ifx\theasciiabstract\relax
\immediate\write\gtoutfile{\theabstract}\else
\immediate\write\gtoutfile{\theasciiabstract}\fi
\immediate\write\gtoutfile{}
\immediate\write\gtoutfile{\noexpand\\}
\immediate\write\gtoutfile{}
\immediate\closeout\gtoutfile}}  

\def\maketitlepage{\makeagttitle\makeheadfile}
\let\makeshorttitle\maketitlepage
\let\maketitle\maketitlepage

\lognumber{20}
\volumenumber{4}
\volumeyear{2004}
\papernumber{20}
\published{10 June 2004}
\pagenumbers{399}{406}
\received{3 March 2004}
\accepted{28 March 2004}

\usepackage{amsmath,amssymb}
\usepackage{graphicx}

\newtheorem{theorem}{Theorem}

\newtheorem{prop}{Proposition}

\newtheorem{cor}{Corollary}

\newcommand{\on}{\operatorname}

\newcommand{\h}{\widehat}

\newcommand{\Spinc}{\on{Spin}^c}

\newcommand{\Z}{\mathbb{Z}}

\newcommand\goth[1]{\mathfrak{#1}}
\newcommand{\s}{\goth{s}}

\renewcommand{\k}{\goth{k}}

\newcommand{\ov}{\overline}

\begin{document}

\author{Olga Plamenevskaya}
\title{Bounds for the Thurston--Bennequin number\\from Floer homology}
\asciititle{Bounds for the Thurston-Bennequin number from Floer homology}
\address{Department of Mathematics, Harvard University\\Cambridge, MA 02138, 
USA}
\email{olga@math.harvard.edu}

\begin{abstract}
Using a knot concordance invariant from the Heegaard Floer theory 
of Ozsv\'ath and Szab\'o, we obtain new bounds for the
Thurston--Bennequin and rotation numbers of Legendrian knots in
$S^3$.  We also apply these bounds to calculate the knot
concordance invariant for certain knots.
\end{abstract}

\asciiabstract{%
Using a knot concordance invariant from the Heegaard Floer theory 
of Ozsvath and Szabo, we obtain new bounds for the
Thurston-Bennequin and rotation numbers of Legendrian knots in
S^3.  We also apply these bounds to calculate the knot
concordance invariant for certain knots.}


\asciikeywords{Legendrian knot, Thurston-Bennequin number, Heegaard
Floer homology}

\maketitle

\section{Introduction}
Let $K$ be a Legendrian knot of genus $g$ in the standard tight contact structure $\xi_{\rm standard}$
on  $S^3$.
It is well-known
that the Thurston--Bennequin and rotation numbers of $K$
satisfy the Thurston--Bennequin inequality
$$
tb(K)+|r(K)|\leq 2g-1.
$$
Although sharp in some cases (e.g.\ right-handed torus knots), in general
this bound is far from optimal. Better bounds can be obtained using
Kauffman and HOMFLY polynomials \cite{FT}, \cite{Ta}.  The Kauffman
polynomial bounds are easily seen to be sharp for left-handed
torus knots; they also allow one to determine the values of the maximal
Thurston--Bennequin number for all two-bridge knots \cite{Ng1}.

In this paper we  use the Ozsv\'ath--Szab\'o knot concordance
invariant $\tau(K)$ introduced in  \cite{tauOS}, \cite{Ra} to
establish a new bound for the Thurston--Bennequin and the rotation number
of a Legendrian knot.
We have

\begin{theorem}\label{tau} For a Legendrian knot $K$ in $(S^3, \xi_{\rm standard})$
 $$
tb(K)+|r(K)|\leq 2\tau(K)-1.
$$
\end{theorem}

For a large class of knots (``perfect'' knots \cite{Ra}),
$\tau(K)=-\sigma(K)/2$,
where $\sigma(K)$ is the signature of the knot (with the sign conventions such that
the right-handed trefoil has signature $-2$).
All alternating knots are perfect \cite{AltOS}, which gives
\begin{cor}\label{alt}
If $K\subset (S^3, \xi_{\rm standard})$ is an alternating Legendrian knot,
then
$$
tb(K)+|r(K)|\leq -\sigma(K)-1.
$$
In particular, for  alternating knots with $\sigma(K)>0$, the
Thurston--Bennequin inequality is not sharp, and $tb(K)$ can never
be positive.
\end{cor}
Note that this bound is usually not sharp even for two-bridge knots
and knots with few crossings
(as can be seen from the calculations in \cite{Ng1}).

It is shown in \cite{tauOS} that $|\tau(K)|\leq g^*(K)$, where
$g^*(K)$ is the four-ball genus of $K$. We therefore recover
a bound due to Rudolph \cite{Ru}:
\begin{cor}\label{g*}
$tb(K)+|r(K)|\leq 2 g^*(K)-1$.
\end{cor}

We prove Theorem \ref{tau} by examining the Heegaard Floer invariants of contact
manifolds obtained by Legendrian surgery.
The Heegaard Floer contact invariants were introduced by Ozsv\'ath and
Szab\'o in \cite{ContOS}; to an oriented contact 3-manifold $(Y,\xi)$
with a co-oriented contact structure $\xi$ they
associate an element $c(\xi)$ of the Heegaard Floer homology group
$\widehat{HF}(-Y)$, defined up to sign. Conjecturally,
the Heegaard Floer contact invariants
are the same as the Seiberg-Witten invariants of contact structures
constructed in \cite{KM}. The definition of $c(\xi)$ uses
an open book decomposition of the contact manifold;
the reader is referred to \cite{ContOS} for the details.

\medskip
{\bf Acknowledgements}\qua I am grateful to Peter Kronheimer and Jake Rasmussen for
illuminating discussions.
\section{The Invariant $\tau(K)$ and Surgery Cobordisms}

In this section we collect the relevant results of Ozsv\'ath, Szab\'o,
and Rasmussen.

For a knot $K\subset S^3$, the invariant $\tau(K)$ is defined via the
Floer complex of the knot; we will
need its interpretation in terms of
surgery cobordisms \cite{tauOS}.

We use notation of \cite{3OS}. Consider the Heegaard Floer group
$\h{HF}(Y)$ of a 3-manifold $Y$, and recall the decomposition
$\widehat{HF}(Y)=\bigoplus_{\s\in\Spinc(Y)}\widehat{HF}(Y,\s)$.
As described in \cite{4OS}, a cobordism $W$ from $Y_1$ to $Y_2$
induces a  map on Floer homology. More precisely, a $\Spinc$
cobordism $(W, \s)$ gives a map
$$
\widehat{F}_{W, \s}:\widehat{HF}(Y_1,\s|Y_1) \to
\widehat{HF}(Y_2,\s|Y_2).
$$

For a knot $K$ in $S^3$ and $n>0$, let $S^3_{-n}(K)$ be obtained by $-n$-surgery on $K$,
and denote by $W$ the cobordism given by the two-handle attachment. The
$\Spinc$ structures on $W$ can be identified with the integers as follows.
Let $\Sigma$ be a Seifert surface for $K$; capping it off inside the attached
two-handle, we obtain a closed surface $\h{\Sigma}$ in $W$.
Let $\s_m$ be the $\Spinc$-structure on $W$ with
$\langle c_1(\s_m), [\h{\Sigma}] \rangle -n=2m$. Accordingly,
the $\Spinc$ structures on  $S^3_{-n}(K)$ are numbered by $[m]\in \Z/n\Z$.
 The cobordism $(W, \s_m)$ induces a map
from  $\h{HF}(S^3)$ to  $\h{HF}(S^3_{-n}(K), [m])$; it will be convenient
to think of $(W, \s_m)$ as a cobordism from $(-S^3_{-n}(K), [m])$ to
$-S^3$, and consider the associated map
$\h{F}_{n,m}:\h{HF}(-S^3_{-n}(K), [m]) \to\h{HF}(-S^3)$.
Now, suppose that $n$ is very large. By the adjunction inequality,
the map $\h{F}_{n,m}$ vanishes for large $m$; moreover, it turns out that
its behavior is controlled by the knot invariant $\tau(K)$:

\begin{prop}{\rm\cite{tauOS, Ra}}\label{m>t}\qua
For all sufficiently large $n$,  the map $\h{F}_{n,m}$
vanishes when $m>\tau(K)$, and is non-trivial when $m<\tau(K)$.
\end{prop}
Note that for $m=\tau(K)$ the map $\h{F}_{n,m}$ might or might not vanish, depending on the knot $K$.

We'll need two more properties of~$\tau(K)$:
\begin{prop}\label{sum}{\rm\cite{tauOS}}

{\rm1)}\qua If the knot $\ov{K}$ is the mirror image of $K$,
then $\tau(\ov{K}) =-\tau(K)$.

{\rm2)}\qua If $K_1\# K_2$ is the connected sum of two knots $K_1$ and $K_2$, then
$\tau(K_1\#K_2)=\tau(K_1)+\tau(K_2)$.
\end{prop}

\section{Contact Invariants and Legendrian Knots}

In this section we use  properties of
the contact invariants  to prove Theorem~\ref{tau}.

Let the contact manifold $(Y_2, \xi_2)$ be obtained from $(Y_1, \xi_1)$ by Legendrian
surgery, and denote by $W$ the corresponding cobordism. As shown in \cite{LS1},
the induced map $\h{F}_{W}$, obtained by summing over $\Spinc$ structures on $W$,
respects the contact invariants; we shall need a slightly stronger statement
for the case of Legendrian surgery on $S^3$, using the canonical $\Spinc$ structure only.

The canonical $\Spinc$ structure $\k$ on the Legendrian surgery cobordism $W$ from $S^3$ to
$S^3_{-n}(K)$ (or, equivalently, from $-S^3_{-n}(K)$ to $-S^3$)
is induced by the Stein structure and determined by the rotation number of $K$,
\begin{equation} \label{canon}
\langle c_1(\k), [\h{\Sigma}] \rangle=r,
\end{equation}
where $\h{\Sigma}$ is the surface
obtained by closing up the Seifert surface of $K$ in the attached Stein handle \cite{Go}.
Let $\s$ be the induced $\Spinc$ structure on $-S^3_{-n}(K)$; $\s$ is the $\Spinc$ structure
associated to $\xi$, and $c(\xi)\in \h{HF}(-S^3_{-n}(K), \s)$.

\begin{prop}[cf.\ \cite{LS1}]\label{lcob}Let  $(W, \k)$  be a
cobordism from $(S^3, \xi_{\rm standard})$ to $(S^3_{-n}(K), \xi)$
induced by Legendrian surgery on $K$, and let $$\h{F}_{W,
\k}:\h{HF}(-S^3_{-n}(K), \s)\to \h{HF}(-S^3)$$ be the associated
map. Then
$$
\widehat{F}_{W, \k}(c(\xi))=c(\xi_{\rm standard}).
$$
\end{prop}
Since $c(\xi_{\rm standard})$ is a generator of $\Z=\h{HF}(S^3)$, it follows that the map
$\h{F}_{W, \k}$ is non-trivial.

\begin{proof}[Proof of Theorem \ref{tau}]  Since
changing the orientation of the knot changes the sign of
its rotation number, it suffices to prove the inequality
\begin{equation}\label{inq}
tb +r \leq 2\tau(K) -1.
\end{equation}
We may also assume that $tb(K)$ is a large
negative number: we can stabilize
the knot (adding kinks to its front projection) to decrease
the Thurston--Bennequin number
and increase the rotation number while keeping $tb+r$ constant.

Writing  $-n=tb-1$  for the coefficient for Legendrian
surgery and setting $r-n=2m$, by (\ref{canon}) we can identify the map $\h{F}_{W,\k}$,
induced by Legendrian surgery, with $\h{F}_{n,m}$ in the notation of Section 2.
By Proposition \ref{lcob}, this map does not vanish, so Proposition \ref{m>t}
implies that $m\leq \tau(K)$, which means that
$$
tb(K)+r(K) \leq 2 \tau(K)+1.
$$
To convert $+1$ into $-1$, we apply this inequality to the knot $K
\# K$. Recalling that $tb(K_1\#K_2)=tb(K_1)+tb(K_2)+1$ and
$r(K_1\#K_2)=r(K_1)+r(K_2)$ and using additivity of $\tau$, we get
$2tb(K)+2r(K)+1 \leq 4 \tau(K) +1$. Then  $tb(K)+r(K) \leq 2
\tau(K)$, and (\ref{inq}) now follows, since $tb(K)+r(K)$ is
always odd (because the numbers $tb(K)-1=
\h{\Sigma}\cdot\h{\Sigma}$ and $\langle c_1(\k), [\h{\Sigma}]
\rangle=r$ have the same parity).
\end{proof}

\noindent {\bf Example 1}\qua  Let $K$ be a $(p, q)$ torus knot. By
\cite{tauOS}, $\tau(K)=\frac{1}{2}(p-1)(q-1)=g(K)$, so Theorem \ref{tau}
reduces to the Thurston--Bennequin inequality, which is actually sharp
in this case. For a  $(-p,q)$ torus knot $\ov{K}$, $\tau(\ov{K})=-\frac{1}{2}(p-1)(q-1)$,
and Theorem \ref{tau} gives $ tb(\ov{K})+|r(\ov{K})|\leq -pq+p+q-2$. Although
stronger than the Thurston--Bennequin inequality, this bound is unfortunately
not sharp: it follows from the Kauffman and HOMFLY polynomial bounds that
$tb(\ov{K})+|r(\ov{K})|\leq -pq$ (and the latter bound is sharp).

\section{An Application: calculating $\tau(K)$}

In this section we use Theorem \ref{tau} to determine the
invariant $\tau(K)$ for certain knots; essentially,
we just give a different proof for some results of \cite{tauOS} and \cite{Li}.

Indeed, a Legendrian representative of $K$ and Theorem \ref{tau}
allows us to find  a lower bound for $\tau(K)$.
An upper bound is given by the  unknotting number of the knot, since
$|\tau(K)|\leq g^*(K)\leq u(K)$. While $u(K)$ is normally hard to determine,
we only need to
look at the unknotting number for some diagram of $K$
to find an upper bound for $\tau(K)$.

\begin{figure}[ht]
\begin{center}
\includegraphics[scale=0.8]{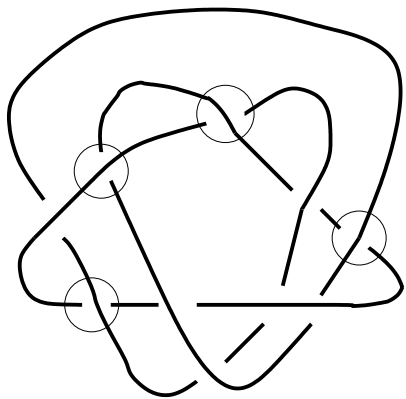} 
\end{center}
\caption{The knot $10_{139}$. Changing the four circled crossings, we obtain an unknot.}
\label{139}
\end{figure}

\begin{figure}[ht]
\begin{center}
\includegraphics[scale=0.9]{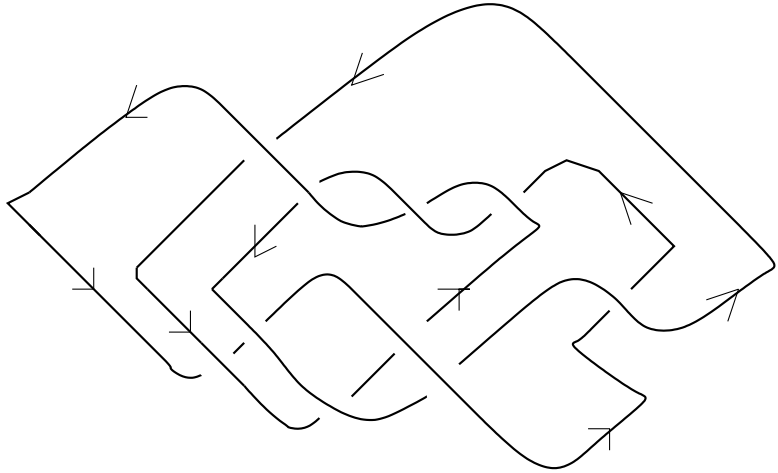} 
\end{center}
\caption{A Legendrian representative of $10_{139}$}
\label{l139}
\end{figure}

\medskip
\noindent {\bf Example 2}\qua \cite{tauOS, Li}\qua
We determine $\tau(K)$ for the knot $K=10_{139}$,
shown on Fig.\ \ref{139}. Changing the four crossings circled on the diagram,
we obtain an unknot. Therefore, $\tau(K)\leq u(K)\leq 4$, so $\tau(K)\leq 4$.
For a lower bound, look at the front projection of the
(oriented) Legendrian representative of $10_{139}$
on Fig.\ \ref{l139}.
The Thurston--Bennequin and the rotation number can be
easily found, since
\begin{eqnarray*}
tb(K)&=&\mbox{ writhe}(K)- \#\mbox{(right cusps)}, \\
r(K)&=&\# \mbox{(upward right cusps)} - \# \mbox{(downward left cusps)}
\end{eqnarray*}
for an oriented front projection. We compute $tb=6$, $r=1$, so $2\tau(K)-1\geq 7$,
and $\tau(K)\geq 4$. It follows that $\tau(K)=4$.

\begin{figure}[h!!t]
\begin{center}
\includegraphics[scale=0.88]{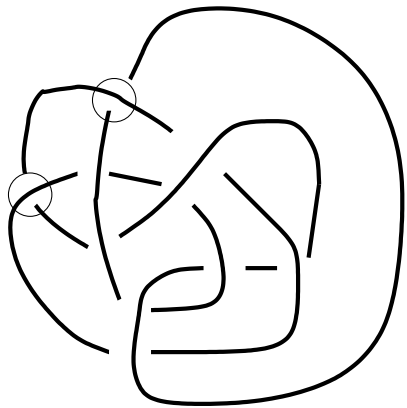} 
\end{center}
\caption{The knot $-10_{145}$. Changing the two circled crossings, we obtain an unknot.}
\label{145}
\end{figure}
\begin{figure}[h!t]
\begin{center}
\includegraphics[scale=0.9]{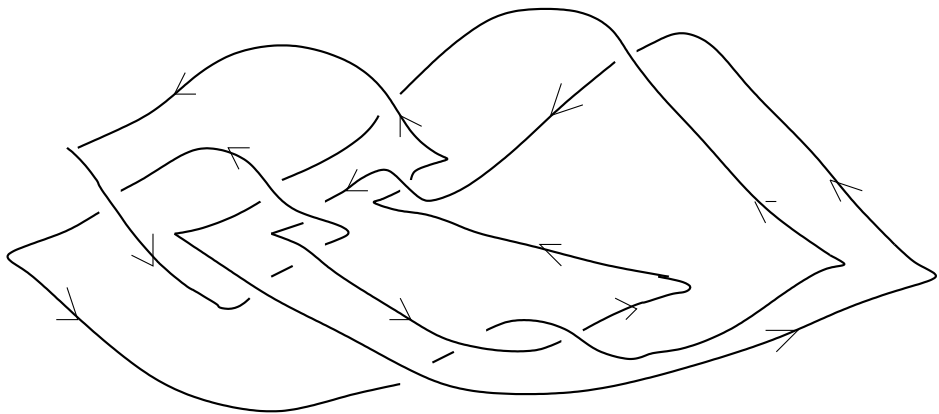} 
\end{center}
\caption{A Legendrian representative of $-10_{145}$}
\label{l145}
\end{figure}

\medskip
\noindent {\bf Example 3}\qua \cite{Li}\qua Using the same idea, we find
 $\tau(K)$ for the knot $K=-10_{145}$,
shown on Fig.\ \ref{145}.
This knot can be unknotted by changing
the two circled crossings, so $\tau(K)\leq u(K)\leq 2$. On the
other hand, for the Legendrian representative shown on Fig.\ \ref{l145},
we compute $tb(K)=2$, $r(K)=1$. Accordingly, $2\tau(K)-1\geq 3$;
it follows that $\tau(K)=2$.

\begin{theorem}{\rm\cite{Li}}\qua Let K be a knot which admits a
Legendrian representive with positive Thurston--Bennequin number,
and let $K_n$ be its $n$-th iterated untwisted positive Whitehead
double. Then $\tau(K_n)=1$.
\end{theorem}
We recall that  an  untwisted positive Whitehead
double for a knot $K\subset S^3$ is constructed by connecting
the knot $K$ and its $0$-push-off $K'$ with a cusp;
here the $0$-push-off is meant to be a copy of $K$, pushed off
in the direction normal
to a Seifert surface for $K$.

\begin{proof}  Clearly, for the Whitehead double of
any knot, we can obtain an unknot by changing one
of the two crossings in the cusp connecting the two copies of the knot.
Then the unknotting number for a Whitehead double cannot be greater than one.
 Now, by a theorem of Akbulut and
Matveyev \cite{AM} the knot $K_n$ has a Legendrian representative
$L_n$ with $tb(L_n)=1$ provided that the original knot $K$
has a Legendrian representive with the  positive Thurston-Bennequin number.
Since $tb(L_n)+|r(L_n)|\leq 2\tau(K_n)-1$, and $\tau(K_n)\leq u(K)\leq 1$,
it follows that $\tau(K_n)=1$.
\end{proof}


\Addresses\recd

\end{document}